\documentstyle[12pt,leqno]{article}
\title{\ }
\author{\ }

\newcommand{\beqa}{\begin{eqnarray}}
\newcommand{\beqan}{\begin{eqnarray*}}
\newcommand{\eeqa}{\end{eqnarray}}
\newcommand{\eeqan}{\end{eqnarray*}}

\newtheorem{theorem}{Theorem}[section]

\newtheorem{lemma}[theorem]{Lemma}

\newtheorem{proposition}{Proposition}[section]

\newtheorem{remark}[theorem]{Remark}

\def\C{{\bf \hbox{\sc I\hskip -7pt C}}} 

\def\Lim{\displaystyle\lim}
\def\Sup{\displaystyle\sup}
\def\R{{\bf \hbox{\sc I\hskip -2pt R}}} 

\def\U{{\bf \hbox{\sc I\hskip -2pt U}}}

\def\N{{\bf \hbox{\sc I\hskip -2pt N}}} 
\def\Z{{\bf Z}} 
\def\Int{\displaystyle\int}
\def\Frac{\displaystyle\frac}
\begin{document}

{\begin{center} {\ } \vskip 2.1cm {\large\bf On the Stability of  swelling porous elastic soils with a single internal
fractional damping*
\let\thefootnote\relax\footnote{*Corresponding author A. BENAISSA, Professor, Ph. D., E-mail: benaissa$_{-}$abbes@yahoo.com}
\\ } \
\\ \ \\
\begin{tabular}{c}
\sc Said Rafa and Abbes Benaissa,
\\ \\
\small  Laboratory of Analysis and Control of PDEs, \qquad \\
\small Djillali Liabes University,\qquad \\
\small P. O. Box 89, Sidi Bel Abbes 22000, ALGERIA.\qquad \\
\end{tabular}
\end{center}
\
\\

\begin{abstract}
{\small
In this paper,we study the well-posedness and polynomial stability to the one-dimensional system in the linear
isothermal theory of swelling porous elastic soils with internal damping of fractional derivative type
We prove well-posedness by using the semigroup theory. Also we establish an optimal decay result by frequency domain
method and Borichev–Tomilov theorem.}
\end{abstract}


{\it ${\cal AMS}$ Classification:} 35Q53, 35Q55, 47J353, 35B35.
\begin{sloppypar}

\renewcommand{\theequation}{\thesection.\arabic{equation}}
\setcounter{equation}{0}

\section{Introduction}
\label{Introduction}
In this paper, we consider the swelling problem in porous elastic soils with fluid saturation give by
\begin{equation}
\left\{
\matrix{
\rho_{z}z_{tt}-a_{1} z_{xx}-a_{2} u_{xx}+\gamma \partial_{t}^{\alpha, \kappa}z(x,t)=0\hfill & \hbox{ in } (0, L)\times(0, +\infty), \hfill \cr
\rho_{u}u_{tt}-a_{3}u_{xx}-a_{2}z_{xx}  = 0\hfill & \hbox{ in } (0, L)\times(0, +\infty),\hfill \cr
z(x,\,0) =  z_{0}(x),\quad z_{t}(x,\,0) = z_{1}(x)\hfill & \hbox{ on } (0, L),\hfill \cr
u(x,\,0) =  u_{0}(x),\quad u_{t}(x,\,0) = u_{1}(x)\hfill & \hbox{ on } (0, L),\hfill \cr
z(0,\,t) =z(L,\,t)=u(0,\,t)=u(L,\,t)= 0,  \hfill & \hbox{ on } (0, +\infty),\hfill \cr}
\right.
\label{eq2}
\end{equation}
where $\gamma> 0$. The notation $\partial_{t}^{\alpha, \kappa}$ stands for the generalized Caputo's fractional derivative
of order $\alpha$, $0<\alpha< 1$, with respect to the time variable. It is defined as follows
$$
\partial_{t}^{\alpha, \kappa}f(t)=\Frac{1}{\Gamma(1-\alpha)}\Int_{0}^{t} (t-s)^{-\alpha}e^{-\kappa(t-s)}\Frac{df}{ds}(s)\, ds, \quad \kappa\geq 0.
$$
A simple model describing the swelling porous elastic soils with fluid  saturation, which was developed in
{\bf\cite{iesa}} and simplified in {\bf\cite{quint}}, is given by a system of coupled hyperbolic equations of the form
$$
\left\{\matrix{\rho_{z}z_{tt}=P_{1x}-G_1+F_1 \hfill &\cr
\rho_{u}u_{tt}=P_{2x}-G_2+F_2\hfill &\cr}\right.
\leqno{(*)}
$$
where the constituents $z$ and $u$ represent the displacement of the fluid and the elastic solid
material, respectively. The positive constant coefficients $\rho_z$ and $\rho_u$ are the densities of each
constituent. The functions $(P_1, G_1, F_1)$ represent the partial tension, internal body forces,
and eternal forces acting on the displacement, respectively. Similar definition holds for
$(P2, G2, F2)$ but acting on the elastic solid. In addition, the constitutive equations of partial
tensions are given by
$$
\pmatrix{P_1 \cr P_2\cr}=\underbrace{\pmatrix{a_1& a_2\cr a_2 & a_3\cr}}_{A}\pmatrix{z_x \cr u_x\cr},
$$
where $a_1, a_3$ are positive constants and $a_2\not= 0$ is a real number. The matrix $A$ is positive
definite in the sense that $a_1a_3\geq a_2^2$. Quintanilla {\bf\cite{quint}} investigated $(1)$ by taking
$$
G_1=G_2=\varsigma(z_t-u_t),\ F_1=a_3z_{xxt},\ F_2=0,
$$
where $\varsigma$ is a positive coefficient, with initial and homogeneous Dirichlet boundary conditions
and obtained an exponential stability result. Similarly, Wang and Guo {\bf\cite{guowa}} considered
$(*)$ with initial and some mixed boundary conditions, taking
$$
G_1=G_2=0,\ F_1=\rho_z \gamma(x) z_{t},\ F_2=0,
$$
where $\gamma(x)$ is an internal viscous damping function with a positive mean. They used the
spectral method to establish an exponential stability result.
Recently, Apalara {\bf\cite{apala}} considered the following
\begin{equation}
\left\{
\matrix{
\rho_{z}z_{tt}-a_{1} z_{xx}-a_{2} u_{xx}=0\hfill & \hbox{ in } (0, 1)\times(0, +\infty), \hfill \cr
\rho_{u}u_{tt}-a_{3}u_{xx}-a_{2}z_{xx}+\Int_{0}^{t}g(t-s)u_{xx}(x, s)\, ds =0\hfill & \hbox{ in } (0, 1)\times(0, +\infty), \hfill \cr
z(x,\,0) =  z_{0}(x),\quad z_{t}(x,\,0) = z_{1}(x)\hfill & \hbox{ on } (0, 1), \hfill \cr
u(x,\,0) =  u_{0}(x),\quad u_{t}(x,\,0) = u_{1}(x)\hfill & \hbox{ on } (0, 1), \hfill \cr
z(0,\,t) =z(1,\,t)=u(0,\,t)=u(1,\,t)= 0\hfill & \hbox{ on } (0, 1)\times(0, +\infty), \hfill \cr}
\right.
\label{eqff}
\end{equation}
where the relaxation function satisfies the condition
$$
g'(t)\leq -\theta(t) g(t),\ t\geq 0
$$
and established a general decay result.

Recently, in {\bf\cite{amma}}, Ammari et al., studied the wave equation with
internal fractional damping. The system considered is as
follows:
\begin{equation}
\left\{
\begin{array}{ll}
u_{tt}(x,t)-\Delta u(x,t)+ \gamma \partial_{t}^{\alpha, \kappa}u(x,t)=0 & \hbox{ in } \Omega\times (0, +\infty),\\
u(x, t)=0  &\hbox{ on } \partial\Omega\times (0, +\infty), \\
u(x, 0)=u_0(x),\quad   u_t(x, 0)=u_1(x) & \hbox{ on }\Omega.
\end{array}
\right.
\nonumber
\end{equation}
The authors proved well-posedness and polynomial asymptotic stability as $t^{-2/(1-\alpha)}$ (for $\kappa> 0$).
Very recently, In {\bf\cite{olive}} (see also {\bf\cite{menno}}), Oleivera et al. considered the porous-elastic system
with two internal fractional dampings defined by
$$
\left\{\matrix{
\rho u_{tt}(x,t)- u_{xx}(x,t)-b\phi_x+\partial_{t}^{\alpha, \kappa}u(x, t)=0\hfill & \hbox{ in } (0, L)\times(0,+\infty), \hfill &\cr
J \phi_{tt}(x,t)- \delta \phi_{xx}(x,t)+bu_x+\eta\phi+\partial_{t}^{\alpha, \kappa}\phi(x, t)=0\hfill & \hbox{ in } (0, L)\times(0,+\infty), \hfill &\cr
u(0, t)=u(L, t)=\phi(0, t)=\phi(L, t)=0\hfill & \hbox{ on } (0, +\infty), \hfill &\cr
u(x,0)=u_0(x), u_{t}(x,0)=u_1(x)\hfill & \hbox{ on } (0, L) \hfill &\cr
\phi(x,0)=\phi_0(x), \phi_{t}(x,0)=\phi_1(x)\hfill & \hbox{ on } (0, L). \hfill &\cr}\right.
\leqno{(PEF)}
$$
They proved, under the condition $\kappa > 0$, that the
energy of system $(PEF)$ decays as time goes to infinity as $t^{-1/(1-\alpha)}$.
To the best of our knowledge, a linear swelling porous elastic soils with one internal
fractional damping has not been studied yet.

This paper is organized as follows. Section 2 briefly outlines the notation and also we reformulate the model
(\ref{eq2}) into an augmented system, coupling the system (\ref{eq2}) with a suitable diffusion equation.
In Section 3, we establish the well-posedness of the system (\ref{eq2}).
In Section 4 we show the lack of exponential stability.
Finally, in Section 5 we obtain the polynomial stability of the corresponding semigroup.

\renewcommand{\theequation}{\thesection.\arabic{equation}}
\setcounter{equation}{0}

\section{Augmented Model}
This section is concerned with the reformulation of the model $(P)$ into an augmented
system. For that, we need the following claims.
\begin{theorem}[\cite{15}]\label{theorem301}
Let $\mu$ be the function
\begin{eqnarray}
\label{601}\mu(y) = |y|^{(2\,\alpha - 1)/2},\quad y\in \R,\quad 0<\alpha<1.
\end{eqnarray}
Then the relation between the {\it input} $U$ and the {\it output}
$O$ if the following system
\begin{eqnarray}
\label{602}&  & \varphi_{t}(y,\,t) + (y^{2}+\kappa)\,\varphi(y,\,t) -
U(t)\,\mu(y) = 0,\quad y\in\R,\quad t>0, \\
\label{603}&  & \varphi(y,\,0) = 0,\\
\label{604}&  & O(t) =
[\pi]^{-1}\,\sin(\alpha\,\pi)\int_{\R}\mu(y)\,\varphi(y,\,t)\
dy
\end{eqnarray}
is given by
\begin{eqnarray}
\label{605}O(t) = I^{1 - \alpha,\,\kappa}U = D^{\alpha,\,\kappa}U,
\end{eqnarray}
where
\begin{eqnarray}
\label{007}[I^{\alpha,\,\kappa}f](t) =
e^{-\,\kappa\,t}\int_{0}^{t}\frac{(t - s)^{\alpha -
1}}{\Gamma(\alpha)}\,e^{\kappa\,s}\,f(s)\ ds.
\end{eqnarray}
\end{theorem}

\begin{lemma}[\cite{acbe}]\label{lemma22}
If $\lambda \in D= \C\backslash ]-\infty, -\kappa]$. Then
\begin{eqnarray*}
\int_{\R}\frac{\mu^{2}(y)}{y^{2} + \kappa + \lambda}\ dy = \frac{\pi}{\sin(\alpha\,\pi)}\,(\kappa + \lambda)^{\alpha - 1}.
\end{eqnarray*}
\end{lemma}
Using the Theorem \ref{theorem301}, we reformulate system (\ref{eq2}), that is, this system can be included
into the augmented model
\begin{equation}
 \left\{
\matrix{
\rho_{z}z_{tt}-a_{1} z_{xx}-a_{2} u_{xx}+ \zeta \int_{-\infty}^{+\infty}\mu(y)\varphi(x,y, t)\, dy =0\hfill & \hbox{ in } (0, L)\times(0, +\infty),\hfill  \cr
\rho_{u}u_{tt}-a_{3}u_{xx}-a_{2}z_{xx}  = 0\hfill & \hbox{ in } (0, L)\times(0, +\infty),  \hfill  \cr
\partial_{t}\varphi(x, y,\,t) + \left(y^{2} + \kappa\right)
\varphi(x, y,\,t) - u_{t}(x,\,t)\,\mu(y) = 0\hfill & \hbox{ in } (0, L)\times \R\times(0, +\infty),  \cr
z(x,\,0) =  z_{0}(x),\quad z_{t}(x,\,0) = z_{1}(x)\hfill & \hbox{ in } (0, L),\hfill  \cr
u(x,\,0) =  u_{0}(x),\quad u_{t}(x,\,0) = u_{1}(x)\hfill & \hbox{ in } (0, L),\hfill  \cr
z(0,\,t) =z(L,\,t)=u(0,\,t)=u(L,\,t)= 0\hfill & \hbox{ in } (0, +\infty),\hfill  \cr
\varphi(x, y, 0)=0\hfill & \hbox{ in } (0, L)\times(-\infty, +\infty),\hfill  \cr}
\right.
\label{aug}
\end{equation}
where we denote $\zeta := \gamma\pi^{-1}\,\sin(\alpha\,\pi)$.
Our first result states that the energy is not conserving.
\begin{lemma}
Let $(z,\,u,\,\varphi)$ be a solution of the system (\ref{aug}). Then, the {\it energy} functional defined by
\begin{eqnarray}
{\cal E}(t) = & \Frac{\rho_{z}}{2}\Int_{0}^{L}|z_{t}|^2dx +
\frac{\rho_{u}}{2}\Int_{0}^{L}|u_{t}|^2dx +
\frac{1}{2}\big(a_{1}-a^{2}_{2}/a_{3}\big)\Int_{0}^{L}|z_{x}|^2dx\qquad \ \ \nonumber \\
\label{e1}& +\frac{1}{2}\Int_{0}^{L}\bigg|\frac{a_{2}}{\sqrt{a_{3}}}z_{x}+\sqrt{a_{3}}u_{x}\bigg|^{2}dx
+ \frac{\zeta}{2}\Int_{0}^{L}\Int_{\R}|\varphi(x, y,\,t)|^{2}\,dy\, dx,
\label{dd20}
\end{eqnarray}
where ${{\cal E}}(t)$ be the energy associated with the coupled system (\ref{aug}) and satisfies
\begin{eqnarray}
\label{e2} \frac{d}{dt}{{\cal E}}(t) = -\
\zeta\Int_{0}^{L}\Int_{\R}\left(y^{2} +\kappa\right)|\varphi(x, y, t)|^{2}\ dy\, dx  \leq 0.
\end{eqnarray}
\end{lemma}
{\bf Proof.}
In fact, multiplying (\ref{eq2})$_{1,\,2}$ by $z_{t}$ and $u_{t}$ respectively, integrating by part over $x\in (0, L)$ and performing straightforward
calculations we obtain
\begin{equation}\label{E4}
\matrix{\Frac{1}{2}\,\frac{d}{dt}\Int_{0}^{L} \left(\rho_{z}|z_{t}|^{2} + a_{1}|z_{x}|^{2} \right)\,dx
+ a_{2}\Re\Int_{0}^{L}\,u_{x}\,{\overline z}_{xt}\,dx\qquad\qquad \hfill \cr
\hfill\qquad\qquad +\zeta\Re\Int_{0}^{L}{\overline z}\Int_{-\infty}^{+\infty}\mu(y)\varphi(x, y, t)\, dy\, dx=0 \cr}
\end{equation}
and
\begin{eqnarray}\label{E5}
\frac{1}{2}\,\frac{d}{dt}\Int_{0}^{L} \left(\rho_{u}|u_{t}|^{2} + a_{3}|u_{x}|^{2} \right)\,dx
+a_{2}\Re\Int_{0}^{L}\,z_{x}{\overline u}_{xt}\,dx = 0.
\end{eqnarray}
Multiplying the third equation in $(P)$ by $\zeta \overline{\varphi}$, integrating over $(0, L)\times(-\infty, +\infty)$ we obtain:
\begin{equation}
\matrix{
 \Frac{\zeta}{2}\Frac{d}{dt}\|\varphi\|_{ L^{2}((0, L) \times(-\infty, +\infty))}^{2}
 + \zeta\Int_{0}^L\int_{-\infty}^{+\infty}(y^{2}+\kappa)|\varphi(x,y,t)|^2 \,dy dx \qquad\qquad\hfill &\cr
- \zeta  \Re \Int_{0}^{L} u_t(x, t) \int_{-\infty}^{+\infty}\mu(y)\overline{\varphi}(x,y,t)dy dx=0.&\cr}
\label{e13}
\end{equation}
From (\ref{dd20}), (\ref{E4}) and (\ref{E5}) we get
$$
\frac{d}{dt}{{\cal E}}(t) = -\
\zeta\Int_{0}^{L}\int_{\R}\left(y^{2} +\kappa\right)|\varphi(x, y, t)|^{2}\ dy\, dx  \leq 0.
$$
\begin{remark}
We have seen that $\frac{d{{\cal E}}(t)}{dt} \neq 0,$ it follows from (\ref{e2}) that it is not energy conserving.
\end{remark}

\renewcommand{\theequation}{\thesection.\arabic{equation}}
\setcounter{equation}{0}
\section{Global existence}
In this section we use results of the semigroup theory of linear operators to obtain an existence theorem of the
system (\ref{aug}). We will use the following standard $L^{2}(0, L)$ space, we are the scalar product and the norm
are denoted by
$$
\langle h,\,g\rangle_{L^{2}(0, L)} = \Int_{0}^{L}h\,\overline{g} \ dx,\quad \|g\|_{L^{2}(0, L)}^{2} = \Int_{0}^{L}|g|^{2}\ dx.
$$
In a similar way, let $L^{2}(\R)$ be the Hilbert space of
all measurable square integrable functions on the real line with the
inner product
\begin{eqnarray*}
\langle h,\,g\rangle_{L^{2}(\R)} = \int_{\R}h\,{\overline g} \ dy, \qquad h,\; g\in L^{2}(\R).
\end{eqnarray*}
We define
\begin{eqnarray*}
H_{0}^{1}(0, L) = \{f\in H^{1}(0, L):\ f(0)=f(L)=0\}.
\end{eqnarray*}
Then
\begin{eqnarray}\label{401}
{\cal H}= H_{0}^{1}(0, L)\times L^{2}(0, L)\times H_{0}^{1}(0, L)\times L^{2}(0, L)\times L^{2}(\R),
\end{eqnarray}
equipped with the inner product given by
\begin{equation}
\matrix{\langle U,\,\widetilde{U}\rangle_{{{\cal H}}} & = &
\rho_z \Int_{0}^{L}w\,\overline{\widetilde{w}}\ dx
+ \rho_u \Int_{0}^{L}v\,\overline{\widetilde{v}}\ dx
+ \left( a_{1}-a^{2}_{2}/a_{3} \right) \Int_{0}^{L}z_{x}\,\overline{\widetilde{z}}_{x}\ dx
 \hfill \cr
&  & + \Int_{0}^{L} \left( \frac{a_{2}}{\sqrt{a_{3}}}z_{x}+\sqrt{a_{3}}u_{x} \right) \overline{\left( \frac{a_{2}}{\sqrt{a_{3}}}\widetilde{z}_x+\sqrt{a_{3}}\widetilde{u}_x \right)}\ dx
  \hfill \cr
&  & + \zeta \Int_{0}^{L}\int_{\R}\varphi\,\overline{\widetilde{\varphi}}\ dy\, dx, \hfill \cr}
\label{402}
\end{equation}
where $ U=(z,\,w,\,u,\,v,\,\varphi)^{T}\ $ and $\ \widetilde{U}
=(\widetilde{z},\,\widetilde{w},\,\widetilde{u},\,\widetilde{v},\,\widetilde{\varphi})^{T}.$
The norm is given by
\begin{eqnarray*}
\|U\|_{{\cal H}}^{2} & = &
\rho_z \int_0^L |w|^2\,dx + \rho_u \int_0^L |v|^2\,dx + \left( a_{1}-a^{2}_{2}/a_{3} \right) \int_0^L |z_x|^2\,dx \\
&  & + \int_0^L \left| \frac{a_{2}}{\sqrt{a_{3}}}z_{x}+\sqrt{a_{3}}u_{x} \right|^2\,dx
+ \zeta \Int_{0}^{L}\int_{\R} |\varphi|^2\ dy\, dx.
\end{eqnarray*}

We now wish to transform the initial boundary value problem (\ref{aug}) to an abstract problem in the Hilbert space ${\cal H}.$ We introduce the functions $z_{t} = w$, $u_{t} = v$ and rewrite the system (\ref{aug}) as the following initial value problem
\begin{eqnarray}\label{403}
\left\{\matrix{
\displaystyle\frac{d}{dt}U(t) = {{\cal A}}U(t), \hfill & \cr
\displaystyle U(0) = U_{0},\quad
\forall\; t > 0,\hfill & \cr}
\right.
\end{eqnarray}
where $U=(z,\,w,\,u,\,v,\,\varphi)^{T}\ $ and $\ U_{0}=(z_{0},\,w_{1},\,u_{0},\,v_{1},\,\varphi_{0})^{T}$.
The operator $\,{\cal A}:{\cal D}({\cal A})\subset {\cal H}\rightarrow {\cal H}$ is given by
\begin{equation}\label{404}
{\cal A}\left(
\begin{array}{c}
z \\
w \\
u \\
v \\
\varphi
\end{array}
\right) =\left(
\begin{array}{c}
w\\
\rho_{z}^{-1} \left(a_1 z_{xx} + a_2 u_{xx}-\zeta \int_{-\infty}^{+\infty}\mu(y)\varphi(x,y)dy\right) \\
v \\
\rho_{u}^{-1} \left(a_3 u_{xx} + a_2 z_{xx} \right) \\
-\,(y^{2} + \kappa)\,\varphi + w(x)\,\mu(y)
\end{array}
\right)
\end{equation}
with domain
\begin{eqnarray*}
{\cal D}({{\cal A}}) =
\bigg\{ (z,\,w,\,u,\,v,\,\varphi)^{T}\in {{\cal H}};\ z,\,u\in H^{2}(0, L),\ w,\,v\in H_{0}^{1}(0, L), \\
 -\,(y^{2} + \kappa)\,\varphi + w(x)\,\mu(y) \in L^{2}({\R}),\, |y|\varphi\in L^{2}({\R})\bigg\}.\\
\end{eqnarray*}

The next result, we will prove the existence and uniqueness of the system solution (\ref{404}).
\begin{proposition}
The operator ${{\cal A}}$ is the infinitesimal generator of a contraction semigroup $\{{{\cal S}}_{{{\cal A}}}(t)\}_{t\geq 0}.$
\end{proposition}

{\bf Proof.}
We will show that ${{\cal A}}$ is a dissipative operator and that the operator $\lambda I - {\cal A}$ is surjective
for any $\lambda > 0$. Then our conclusion follows using the well known the Lumer-Phillips
Theorem. We observe that if $U = (z,\,w,\,u,\,v,\,\varphi)^{T}\in {{\cal D}}({{\cal A}}),$ then using
(\ref{e2}) and (\ref{403}), we get
\begin{eqnarray}\label{406}
\Re\langle{{\cal A}}U,\,U\rangle_{{\cal H}} =
-\zeta\Int_{0}^{L}\Int_{\R}\left(y^{2} + \kappa\right)|\varphi(x, y)|^{2}\, dy\, dx
\leq 0.
\end{eqnarray}
In fact, for any $U\in D({\cal A})$ we known that
$$
{\cal E}(t)=\Frac{1}{2}\|U\|_{{\cal H}}^2
$$
then we have (\ref{406}). It follows that ${{\cal A}}$ is a dissipative operator.

Next, we will prove that the operator $\lambda\,I - {{\cal A}}$ is surjective for $\lambda >0.$ For this purpose,
let $F = (f_{1},\,f_{2},\,f_{3},\,f_{4},\,f_{5})^{T}\in {{\cal H}}$, we seek
$U=(z,\,w,\,u,\,v,\,\varphi)^{T}\in {{\cal D}}({{\cal A}})$ such that $(\lambda\,I - {{\cal A}})U=F,$ that is,
\begin{eqnarray}
\label{407} \left\lbrace
\begin{array}{l}
\lambda\,z - w = f_{1}, \\
\lambda\,\rho_z\,w - a_1 z_{xx} - a_2 u_{xx}+\zeta \int_{-\infty}^{+\infty}\mu(y)\varphi(x,y)dy = \rho_z\,f_{2}, \\
\lambda\,u - v = f_{3}, \\
\lambda\,\rho_u\,z - a_3 u_{xx} - a_2 z_{xx} = \rho_u\,f_{4}, \\
\lambda\,\varphi + (y^{2} + \kappa)\,\varphi - w(x)\,\mu(y) = f_{5}.
\end{array}
\right.
\end{eqnarray}
From (\ref{407})$_{5}$ we obtain
\begin{eqnarray}
\label{408}
\begin{array}{l}
\varphi(y) = \frac{f_{5}(x, y) \ +\ w(x)\,\mu(y)}{y^{2}\ +\ \kappa\ +\ \lambda}
\end{array}
\end{eqnarray}
and by (\ref{407})$_{1,\,3}$ we have
\begin{eqnarray}
\label{409}\left\lbrace
\begin{array}{l}
w = \lambda\,z - f_{1}, \\
v = \lambda\,u - f_{3}.
\end{array}
\right.
\end{eqnarray}
On the other hand, replacing $(\ref{409})_{1,\,2}$ into (\ref{407})$_{2,\,4}$ and (\ref{408}) into $(\ref{407})_{2}$
respectively, yields
\begin{equation}
\label{410}
\left\{
\matrix{
\lambda^{2}\,\rho_z\,z - a_1\,z_{xx} - a_2\,u_{xx}+\gamma\lambda(\lambda+\kappa)^{\alpha-1}z
= \rho_z\,f_{2} + \lambda\,\rho_z\,f_{1}\hfill &\cr
+\gamma (\lambda+\kappa)^{\alpha-1}f_1-\zeta \int_{-\infty}^{+\infty}\frac{\mu(y)f_{5}(x, y)}{y^2+\lambda+\kappa}dy,
\hfill &\cr
\lambda^{2}\,\rho_u\,v - a_3\,u_{xx} - a_2\,z_{xx} = \rho_u\,f_{4} + \lambda\,\rho_u\,f_{3}.\hfill &\cr}\right.
\end{equation}
To solve (\ref{410}) is equivalent to finding $z,\,u\in H^{2}(0, L)\cap H_{0}^{1}(0, L)$ such that
\begin{eqnarray}
\label{411}\left\lbrace
\begin{array}{l}
\displaystyle{\int_{0}^{L}}\left( (\rho_z\lambda^{2}+\gamma\lambda(\lambda+\kappa)^{\alpha-1})\,\,z - a_1\,z_{xx}
- a_2\,u_{xx} \right) {\overline{\tilde{z}}}\ dx\\
= \displaystyle{\int_{0}^{L}}\left( \rho_z\,f_{2} + (\rho_z\lambda+\gamma (\lambda+\kappa)^{\alpha-1})\,\,f_{1} \right)
{\overline{\tilde{z}}}\ dx,
-\zeta \int_{0}^{L}{\overline{\tilde{z}}}\int_{-\infty}^{+\infty}\frac{\mu(y)f_{5}(x, y)}{y^2+\lambda+\kappa}dy\, dx, \\
\displaystyle{\int_{0}^{L}}\left( \lambda^{2}\,\rho_u\,v - a_3\,u_{xx} - a_2\,z_{xx} \right) {\overline{\tilde{u}}}\ dx
= \displaystyle{\int_{0}^{L}}\left( \rho_u\,f_{4} + \lambda\,\rho_u\,f_{3} \right) {\overline{\tilde{u}}}\ dx,
\end{array}
\right.
\end{eqnarray}
for all $\tilde{z},\ \tilde{u} \in H_{0}^{1}(0, L)$. We consider (\ref{411})$_{1}$, then integrating by parts,
using (\ref{aug})$_{8}$ and from (\ref{408})
it follows that
\begin{equation}\label{412}
\matrix{(\rho_z\lambda^{2}+\gamma\lambda(\lambda+\kappa)^{\alpha-1}) \Int_{0}^{L} z\ {\overline{\tilde{z}}}\ dx
 + a_1 \Int_{0}^{L}z_{x}\,\tilde{z}_{x}\ dx+ a_2 \Int_{0}^{L}u_x\,{\overline{\tilde{z}}}_x\ dx\hfill &\cr
= \int_{0}^{L}\left( \rho_z\,f_{2} + (\rho_z\lambda+\gamma (\lambda+\kappa)^{\alpha-1})\,\,f_{1} \right) {\overline{\tilde{z}}}\ dx,
-\zeta \int_{0}^{L}{\overline{\tilde{z}}}\int_{-\infty}^{+\infty}\frac{\mu(y)f_{5}(x, y)}{y^2+\lambda+\kappa}dy\, dx.
\hfill &\cr}
\end{equation}
In a similar way we estimate (\ref{411})$_{2}$,
\begin{eqnarray}
\lambda^{2}\,\rho_u \Int_{0}^{L} u\ {\overline{\tilde{u}}}\ dx
+ a_3 \Int_{0}^{L}u_{x}\,{\overline{\tilde{u}}}_{x}\ dx
+ a_2 \Int_{0}^{L}u_x\,{\overline{\tilde{z}}}_x\ dx  \nonumber
\\
 = \Int_{0}^{L} \left(\rho_u\,f_{4} + \lambda\,\rho_u\,f_{3} \right){\overline{\tilde{u}}}\ dx.
 \label{4141}
\end{eqnarray}
The system (\ref{412})-(\ref{4141}) is equivalent to the problem
\begin{eqnarray}\label{415}
a\big(\left(z, u), (\tilde{z}, \tilde{u}\right)\big) = {{\cal L}}(\tilde{z}, \tilde{u}),
\end{eqnarray}
where the sesquilinear form $a: \left[H_{0}^{1}(0, L)\times H_{0}^{1}(0, L)\right]^{2} \rightarrow\R$ and the antilinear
form ${{\cal L}}: H_{0}^{1}(0, L)\times H_{0}^{1}(0, L) \rightarrow\R$ are defined by
$$
\matrix{a\big(\left(z, u), (\tilde{z}, \tilde{u}\right)\big)& =&
\, (\rho_z\lambda^{2}+\gamma\lambda(\lambda+\kappa)^{\alpha-1}) \displaystyle \Int_{0}^{L} z\ {\overline{\tilde{z}}}\, dx
+ \rho_u \lambda^{2}\Int_{0}^{L}u\ {\overline{\tilde{u}}}  dx
+ a_1 \Int_{0}^{L}z_{x}\,{\overline{\tilde{z}}}_{x}\ dx \cr
&&+ a_3 \Int_{0}^{L}u_{x}\,\tilde{u}_{x}\ dx
 + a_2 \Int_{0}^{L} \left( u_x\,{\overline{\tilde{z}}}_x + z_x\,{\overline{\tilde{u}}}_x \right) dx
 \hfill \cr}
$$
and
\begin{eqnarray*}
{{\cal L}}(\tilde{z},\,\tilde{u}) & = &
 \int_{0}^{L}\left( \rho_z\,f_{2} + (\rho_z\lambda+\gamma (\lambda+\kappa)^{\alpha-1})\,\,f_{1} \right) {\overline{\tilde{z}}}\ dx
+\Int_{0}^{L} \left(\rho_u\,f_{4} + \lambda\,\rho_u\,f_{3} \right){\overline{\tilde{u}}}\ dx \nonumber
\\
&& -\zeta \int_{0}^{L}{\overline{\tilde{z}}}\int_{-\infty}^{+\infty}\frac{\mu(y)f_{5}(x, y)}{y^2+\lambda+\kappa}dy\, dx.
\label{417}
\end{eqnarray*}
It is easy to verify to $a$ is continuous and coercive, and ${{\cal L}}$ is continuous. So applying the Lax-Milgram
Theorem, we deduce for all
$(\tilde{z},\,\tilde{u})\in H_{0}^{1}(0, L)\times H_{0}^{1}(0, L)$ the problem (\ref{415}) admits a unique solution
$(z,\,u)\in H_{0}^{1}(0, L) \times H_{0}^{1}(0, L).$ Using classical elliptic regularity, it follows from (\ref{412})-(\ref{4141}) that
$(z,\,u)\in H^{2}(0, L)\times H^{2}(0, L)$. Therefore, the operator $\lambda\,I - {{\cal A}}$ is surjective for any $\lambda > 0.$\\

Now, using the Hille-Yosida theorem we have the following results.
\begin{theorem}
\begin{itemize}
\item[(a)]If $U_{0}\in {{\cal H}},$ then the system (\ref{403}) has a unique mild solution
\begin{eqnarray}\label{418}
U\in C\left([0, +\infty),\ {{\cal H}}\right).
\end{eqnarray}
\item[(b)]If $U_{0}\in {{\cal D}}({{\cal A}}),$ then the system (\ref{403}) has a unique strong solution
\begin{eqnarray}\label{419}
U\in C\left([0, +\infty),\ {{\cal D}}({{\cal A}})\right)\cap C^{1}\left([0, +\infty),\ {{\cal H}}\right).
\end{eqnarray}
\end{itemize}
\end{theorem}

\renewcommand{\theequation}{\thesection.\arabic{equation}}
\setcounter{equation}{0}
\section{The lack of exponential stability}

The methods we use to show the lack of exponential stability can be found in \cite{15}. For this, the following results are fundamental.
\begin{theorem}\cite{15}\label{theorem8.2}
The following two statements are equivalent in a Hilbert space:
\begin{enumerate}
\item The semigroup $\{{{\cal S}}(t):\ t\geq 0\}$ is exponentially stable:
\begin{eqnarray}
\exists\; M>0,\quad \exists\; a>0,\quad \|{{\cal S}}(t)\| \leq M\,e^{-\,a\,t}, \label{801}
\end{eqnarray}
where the norm used is the operator norm induced by the inner product (\ref{402}).
\item The set $\{i\,\lambda:\ \lambda\in\R\}$ belongs to the resolvent set of the operator ${{\cal A}}$, i.e, $i\R \subset \varrho ({\cal A})$, and
\begin{eqnarray}
\sup\left\{\left\|\left(i\,\lambda\,I - {{\cal A}}\right)^{-1}\right\|,\quad \lambda \in\R\right\} < \infty. \label{802}
\end{eqnarray}
\end{enumerate}
\end{theorem}

Equipped with the foregoing theorem, we are now able to state and prove this claim of ours.
\begin{theorem}
\label{theorem8.3}The semigroup $\{{{\cal S}}(t):\ t\geq 0\}$ generated by the operator ${{\cal A}}$ is not exponentially stable.
\end{theorem}
{\bf Proof.}
We devide it into two cases: $\kappa = 0$ and $\kappa \neq 0$.

\textbf{Case 1: $\kappa =0.$}

We will show that $0 = i\,\lambda \not\in \varrho({{\cal A}})$. In fact, we consider
$(\sin(\frac{\pi x}{L}),\,0,\,0,\,0,\,0)^{T}\in {{\cal H}},$ then denoting by $(z,\,w,\,u,\,v,\,\varphi)^{T}$
the image of $(\sin(\frac{\pi x}{L}),\,0,\,0,\,0,\,0)^{T}$ by ${{\cal A}}^{-1},$
we can see that $\varphi(x, y) = -\sin(\frac{\pi x}{L})|y|^{(2\,\alpha - 5)/2}$.
But $\varphi\not\in L^{2}((0, L)\times\R)$ for $\alpha\in (0,\,1)$.
This way $(z,\,w,\,u,\,v,\,\varphi)^{T}\not\in {{\cal D}}({{\cal A}}),$ which in turns implies
that $i\,\lambda = 0$ is not in the resolvent set of ${{\cal A}}.$

\textbf{Case 2: $\kappa \neq 0.$}

Using the Theorem \ref{theorem8.2}, it suffices to prove that $\|U\|$ is unbounded for $\lambda\in\R.$

In fact, let $\lambda$ be an eigenvalue of ${{\cal A}}$ with associated
$U = (z,\,w,\,u,\,v,\,\varphi)^{T}\in{{\cal D}}({{\cal A}})$ then ${{\cal A}}U = \lambda\,U,$ is equivalent
\begin{eqnarray}
\label{803.b}
\left\lbrace
\begin{array}{l}
\lambda\,z - w = 0,
\\
\lambda\,\rho_z\,w - a_1\,z_{xx} - a_2\,u_{xx}+\zeta \int_{-\infty}^{+\infty}\mu(y)\varphi(x,y)dy = 0,
\\
\lambda\,u - v = 0,
\\
\lambda\,\rho_u\,v - a_3\,u_{xx} - a_2\,z_{xx} = 0,
\\
(\lambda + \kappa +y^{2})\,\varphi(y) - w(x)\,\mu(y) = 0,
\\
u(0) = 0,\quad z(0) = 0,\quad u(L)=0,\quad z(L)=0.
\end{array}
\right.
\end{eqnarray}
From $(\ref{803.b})_{1}$-$(\ref{803.b})_{2}$ and $(\ref{803.b})_{3}$-$(\ref{803.b})_{4}$ for such $\lambda$, we find
\begin{eqnarray}
\label{803.c}
\left\{\matrix{\displaystyle\left(\lambda^{2}+\frac{\gamma}{\rho_z}\lambda(\lambda+\kappa)^{\alpha-1}\right)\,z - \frac{a_1}{\rho_z}\,z_{xx} - \frac{a_2}{\rho_z}\,u_{xx} = 0,\hfill & \cr
\displaystyle\lambda^{2}\,u - \frac{a_3}{\rho_u}\,u_{xx} - \frac{a_2}{\rho_u}\,z_{xx} = 0.\hfill & \cr}
\right.
\end{eqnarray}
There exists an increasing sequence $(\mu_n)_{n\geq 1}$
tending to infinity and an orthonormal basis $(e_n)_{n\geq 1}$ of $H_{0}^{1}(0, L)$ such that
\begin{equation}
-\partial_{xx}e_n=\mu_n^2 e_n.
\label{zaa11}
\end{equation}
For every $n\geq 1$, there exists $(B_n,C_n)\not= (0, 0)$ such that the eigenvector $U$ of ${\cal A}$ is of the form
\begin{equation}
z=B_n e_n,\quad w=\lambda B_n e_n,\quad u=C_n e_n,\quad v=\lambda C_n e_n.
\label{zaa12}
\end{equation}
Inserting (\ref{zaa12}) in (\ref{803.c}) and using (\ref{zaa11}), we obtain
\begin{equation}
\left\{\matrix{(\lambda^2+\frac{\gamma}{\rho_z}\lambda(\lambda+\kappa)^{\alpha-1}+\frac{a_1}{\rho_z}\mu_n^2)B_n e_n+
\frac{a_2}{\rho_z}\mu_n^2C_n e_n=0, \hfill &\cr
\frac{a_2}{\rho_u}\mu_n^2 B_n e_n+
(\lambda^2+\frac{a_3}{\rho_u}\mu_n^2)C_n e_n=0 \hfill &\cr}\right.
\label{zaa13}
\end{equation}
which has a non-trivial solution $(B_n,C_n)\not=(0, 0)$ if and only if $\lambda$ is a solution of the equation
\begin{equation}
(\lambda^2+\frac{\gamma}{\rho_z}\lambda(\lambda+\kappa)^{\alpha-1}+\frac{a_1}{\rho_z}\mu_n^2)
(\lambda^2+\frac{a_3}{\rho_u}\mu_n^2)-\frac{a_2^2}{\rho_z\rho_u}\mu_n^4=0.
\label{zaa14}
\end{equation}
Then
\begin{equation}
\lambda^4+(\frac{a_1}{\rho_z}+\frac{a_3}{\rho_u})\mu_n^2\lambda^2+\frac{\gamma}{\rho_z}\lambda^3(\lambda+\kappa)^{\alpha-1}
+\frac{\gamma}{\rho_z}\lambda(\lambda+\kappa)^{\alpha-1}\frac{a_3}{\rho_u}\mu_n^2+(\frac{a_1a_3}{\rho_z\rho_u}
-\frac{a_2^2}{\rho_z\rho_u})\mu_n^4=0,
\label{zaa16}
\end{equation}
that we refer to as the characteristic equation associated with the eigenvalue $\mu_n^2$
of ${-\partial_{xx}}$. The four roots
of this equation are eigenvalues of ${\cal A}$ and called the eigenvalues of ${\cal A}$ corresponding to
$\mu_n$ of ${-\partial_{xx}}$.
We also have the following result.
\begin{lemma}
Let $\lambda_n=\lambda_{j, n}^{\pm}, j = 1, 2$ be one of the fourth eigenvalues of ${\cal A}$ corresponding to $\mu_n$.
Then, there exists two positive constants $m, M$, such that, for $n$ large enough,
\begin{equation}
m\leq \left|\Frac{\lambda_n}{\mu_n}\right|\leq M.
\label{zaa15}
\end{equation}
\end{lemma}
{\bf Proof.}
Set $Z_n =\frac{\lambda_n}{\mu_n}$.
Then, from (\ref{zaa16}), one has that $Z_n$ is one of the four roots of the polynomial
$f_n$ of degree four given by
$$
f_n(z)=z^4+(\frac{a_1}{\rho_z}+\frac{a_3}{\rho_u})z^2+\frac{\gamma}{\rho_z}\mu_n^{\alpha-2}z^3(z+\frac{\kappa}{\mu_n})^{\alpha-1}
+\frac{\gamma}{\rho_z}\mu_n^{\alpha-2}z(z+\frac{\kappa}{\mu_n})^{\alpha-1}\frac{a_3}{\rho_u}+(\frac{a_1a_3}{\rho_z\rho_u}
-\frac{a_2^2}{\rho_z\rho_u}).
$$
Let $g$ be the the polynomial of degree four given by
$$
g(z) = z^4+(\frac{a_1}{\rho_z}+\frac{a_3}{\rho_u})z^2+(\frac{a_1a_3}{\rho_z\rho_u},
-\frac{a_2^2}{\rho_z\rho_u})
\leqno{(EA)}
$$
which has exactly four non zero roots. Since the coefficients of $f_n$ converge to those of $g$ as n tends to infinity, one gets
the result.
\begin{lemma}(Asymptotic behavior of the large eigenvalues of ${\cal A}$)
The large eigenvalues of ${\cal A}$ can be split into two families $(\lambda_n^j)_{n\in \Z, |n|\geq n_0}, j=1, 2$, ($n_0\in \N$ chosen
large enough). The following asymptotic expansions hold:
\begin{equation}
\matrix{\lambda_n^1&=&i l_1\mu_n+\frac{\tilde\alpha_1}{\mu_n^{1-\alpha}}+\frac{\beta_1}{\mu_n^{1-\alpha}}+o\left(\frac{1}{\mu_n^{1-\alpha}}\right),
n\geq N, {\tilde\alpha}_1\in i\R, \beta_1\in \R,\beta_1< 0, n\geq n_0\hfill \cr
\lambda_n^2&=&i l_3\mu_n+\frac{\tilde\alpha_2}{\mu_n^{1-\alpha}}+\frac{\beta_2}{\mu_n^{1-\alpha}}+o\left(\frac{1}{\mu_n^{1-\alpha}}\right),
n\geq N, {\tilde\alpha}_2\in i\R, \beta_2\in \R,\beta< 0, n\geq n_0\hfill \cr
&&\lambda_n^1=\overline{\lambda_{-n}^1}, \ \lambda_n^2=\overline{\lambda_{-n}^2} \hbox{ if } n\leq -N.\hfill \cr}
\label{0299hzz}
\end{equation}
and these two roots are simple.
\label{llle}
\end{lemma}
{\bf Proof.}
The equation $(EA)$ has exactly four non zero roots. That is
\begin{eqnarray}\label{raizes}
\left\{\matrix{
\displaystyle \lambda_{0n}^1= i\sqrt{(\frac{a_{1}}{\rho_{z}}+\frac{a_{3}}{\rho_{u}})
+\sqrt{(\frac{a_{1}}{\rho_{z}}-\frac{a_{3}}{\rho_{u}})^{2}+\frac{4a_{2}^{2}}{\rho_{u}\rho_{z}}}}\mu_n, \hfill & \cr
\displaystyle \lambda_{0n}^2=-\lambda_{0n}^1,\hfill & \cr
\displaystyle \lambda_{0n}^3=i\sqrt{(\frac{a_{1}}{\rho_{z}}+\frac{a_{3}}{\rho_{u}})
-\sqrt{(\frac{a_{1}}{\rho_{z}}-\frac{a_{3}}{\rho_{u}})^{2}+\frac{4a_{2}^{2}}{\rho_{u}\rho_{z}}}}\mu_n,\hfill & \cr
\displaystyle \lambda_{0n}^4=-\lambda_{0n}^3,\hfill & \cr}
\right.
\end{eqnarray}
{\bf Step 3.} From Step 2, we can write
\begin{equation}
\lambda_{jn}=\lambda_{jn}^0+\varepsilon_{jn}, \quad j=1, 3.
\label{030e}
\end{equation}
Substituting (\ref{030e}) into (\ref{zaa16}), using that $f(\lambda_{jn})=0$, we get:
$$
4\varepsilon_{jn}(\lambda_{jn}^{0})^3+2 m \mu_n^2\varepsilon_{jn}\lambda_{jn}^{0}
+\frac{\gamma}{\rho_{z}}(\lambda_{jn}^{0})^{2+\alpha}+\frac{\gamma}{\rho_{z}}\frac{a_3}{\rho_{u}}\mu_n^2(\lambda_{jn}^{0})^{\alpha}=0
$$
where
$$
m=\frac{a_{1}}{\rho_{z}}+\frac{a_{3}}{\rho_{u}}
$$
and hence
\begin{equation}
\matrix{
\varepsilon_{jn}&=&-i^{\alpha-1}\Frac{(l_j^2-\frac{a_3}{\rho_{u}})}{(4l_j^2-2m)}\Frac{\gamma}{\rho_{z}(l\mu_n)^{1-\alpha}}+o\left(\frac{1}{\mu_n^{1-\alpha}}\right)\hfill \cr
&=&-\Frac{(l_j^2-\frac{a_3}{\rho_{u}})}{(4l^2-2m)}\Frac{\gamma}{\rho_{z}(l_j\mu_n)^{1-\alpha}}\left(\cos(1-\alpha)\frac{\pi}{2}-i\sin(1-\alpha)\frac{\pi}{2}\right)
+o\left(\frac{1}{\mu_n^{1-\alpha}}\right)\hbox{ for }n\succeq 0,
\hfill \cr}
\label{055}
\end{equation}
where $l_1=l_{+}$ and $l_3=l_{-}$ and
$$
l_{\pm}=\sqrt{(\frac{a_{1}}{\rho_{z}}+\frac{a_{3}}{\rho_{u}})
\pm\sqrt{(\frac{a_{1}}{\rho_{z}}-\frac{a_{3}}{\rho_{u}})^{2}+\frac{4a_{2}^{2}}{\rho_{u}\rho_{z}}}}.
$$
From (\ref{055}) we have in that case $|\mu_n|^{1-\alpha}\Re\lambda_n\sim \beta$, with
$$
\beta=-\Frac{(l_j^2-\frac{a_3}{\rho_{u}})}{(4l_j^2-2m)}\Frac{\gamma}{\rho_{z}l_j^{1-\alpha}}\cos(1-\alpha)\frac{\pi}{2}.
$$
The operator ${\cal A}$ has a non exponential decaying branch of eigenvalues. Thus the proof is complete.

\hfill$\Box$\\

\renewcommand{\theequation}{\thesection.\arabic{equation}}
\setcounter{equation}{0}
\section{Asymptotic behavior}
\subsection{Strong stability of the system}
In this subsection, we will use a general criteria due to Arendt-Batty and Lyubich-Vu (see \cite{arendt}, \cite{lyub})
to prove the strong stability of the C$_{0}$-semigroup $e^{t{{\cal A}}}$. associated to the system (\ref{eq2})
in the absence of the compactness of the resolvent of ${{\cal A}}$.
The main result is the following theorem:
\begin{theorem}
The C$_{0}$-semigroup $e^{t{{\cal A}}}$ is strongly stable in ${{\cal H}},$ that is, for all $U_{0}\in {{\cal H}},$ the solution of (\ref{403}) satisfies
\begin{eqnarray*}
\lim_{t\rightarrow +\infty} \|e^{t{{\cal A}}}U_{0}\|_{{\cal H}} = 0.
\end{eqnarray*}
\end{theorem}

To prove our first result, we use the following theorem.
\begin{theorem}\cite{arendt}-\cite{lyub} \label{arendt}
Let ${{\cal A}}$ be the generator of a uniformly bounded $C_{0}$-semigroup $\{{{\cal S}}(t)\}_{t\geq 0}$ on a Hilbert space ${\cal H}$. If
\begin{enumerate}
\item[(i)]If $\sigma({{\cal A}}) \cap i\,\R$ is at most a countable set, where $\sigma({{\cal A}})$ denotes the spectrum of ${\cal A}$;
\item[(ii)]If $\sigma_r({{\cal A}}) \cap i\,\R = \emptyset$, where $\sigma_r({{\cal A}})$ denotes the set of residual spectrum of ${\cal A}$.
\end{enumerate}
Then the semigroup $\{{{\cal S}}(t)\}_{t\geq 0}$ is asymptotically stable, that is,
\begin{eqnarray*}
\|{{\cal S}}(t)w\|_{{{\cal H}}}\rightarrow 0\quad\mbox{as} \quad t\rightarrow \infty,
\end{eqnarray*}
for any $w \in {{\cal H}}$.
\end{theorem}
To prove this result, we need some lemmas.
\begin{lemma}
We have
\begin{eqnarray}
\label{506}\sigma({{\cal A}}) \cap i\,\R^* = \sigma({{\cal A}}) \cap \{i\,\lambda,\
\lambda\in\R,\ \lambda\neq 0\} = \emptyset.
\end{eqnarray}
\label{lemma1}
\end{lemma}

{\bf Proof.}
By contradiction. We suppose that there
$\lambda\in\R,$ $\lambda\neq 0$ and $U\neq 0,$ such that
\begin{equation}
{{\cal A}}U = i\,\lambda\,U,
\label{ss1}
\end{equation}
that is, $\left(i\,\lambda - {{\cal A}}\right)U = 0.$ Then
\begin{eqnarray}
\label{507} \left\lbrace
\begin{array}{l}
i\lambda\,z - w = 0, \\
i\lambda\,\rho_z\,w - a_1 z_{xx} - a_2 u_{xx}+\zeta \int_{-\infty}^{+\infty}\mu(y)\varphi(x,y)dy = 0, \\
i\lambda\,u - v = 0, \\
i\lambda\,\rho_u\,v - a_3 u_{xx} - a_2 z_{xx} = 0, \\
i\lambda\,\varphi + (y^{2} + \kappa)\,\varphi - w(x)\,\mu(y) = 0.
\end{array}
\right.
\end{eqnarray}
Next, a straightforward computation gives
\begin{equation}
\Re\langle{\cal A}U, U\rangle =-\zeta\int_{0}^{L}\int_{\R}\left(y^{2} + \kappa\right)|\varphi(x, y)|^{2}\ dy\, dx.
\label{ss2}
\end{equation}
Then, using Eq.(\ref{ss2}) and Eq.(\ref{ss1}) we deduce that
\begin{equation}
\varphi\equiv 0 \hbox{ in }(0, L)\times (-\infty, +\infty).
\label{ss3}
\end{equation}
From Eq.$(\ref{507})_5$ and Eq.$(\ref{507})_1$, we have
\begin{equation}
w\equiv 0 \hbox{ and } z\equiv 0 \hbox{ in }(0, L).
\label{ss4}
\end{equation}
As $a_2\not= 0$, we have
$$
u_{xx}=0\hbox{ in }(0, L).
$$
Then $u\equiv 0$. From Eq.$(\ref{507})_3$, we have $v\equiv 0$.
and consequently, $U = 0$.

\noindent
{\bf Step 2:} $\lambda = 0$. The system Eq.(\ref{507}) becomes
\begin{eqnarray}
\label{5077} \left\lbrace
\begin{array}{l}
- w = 0, \\
- a_1 z_{xx} - a_2 u_{xx}+\zeta \int_{-\infty}^{+\infty}\mu(y)\varphi(x,y)dy = 0, \\
- v = 0, \\
 - a_3 u_{xx} - a_2 z_{xx} = 0, \\
(y^{2} + \kappa)\,\varphi - w(x)\,\mu(y) = 0.
\end{array}
\right.
\end{eqnarray}
Because
\begin{equation}
\varphi\equiv 0 \hbox{ in }(0, L)\times (-\infty, +\infty).
\label{ss30}
\end{equation}
We deduce that
\begin{eqnarray}
\label{5078} \left\lbrace
\begin{array}{l}
- w = 0, \\
- a_1 z_{xx} - a_2 u_{xx}= 0, \\
- v = 0, \\
 - a_3 u_{xx} - a_2 z_{xx} = 0.
\end{array}
\right.
\end{eqnarray}
As $a_2^2\not= a_1a_3$, we have $u_{xx}\equiv 0$ and $z_{xx}\equiv 0$.
Then $u\equiv 0$ and $z\equiv 0$ and consequently, $U = 0$.

In the Theorem \ref{arendt}, the condition $(i)$ holds if we show that any point $\sigma({{\cal A}})\cap \{i\,\R\}$ is at most a countable set.

In fact, we will show that the operator $i\,\lambda\,I - {{\cal A}}$ is surjective for $\lambda\neq 0.$
For this purpose, given $F=(f_{1},\,f_{2},\,f_{3},\,f_{4},\,f_{5})^{T} \in {{\cal H}}$, we seek
$U=(z, w, u,\,v,\varphi)^{T}\in {{\cal D}}({{\cal A}})$ which is solution of
$(i\,\lambda\,I - {{\cal A}})U=F,$ that is, the entries of $U$ satisfy the system of equations
\begin{eqnarray}
\label{512} \left\lbrace
\begin{array}{l}
i\,\lambda\,z - w = f_1, \\
i\,\lambda\,\rho_z\,w - a_1\,z_{xx} - a_2\,u_{xx}+\zeta \int_{-\infty}^{+\infty}\mu(y)\varphi(x,y)dy = \rho_z\,f_2, \\
i\,\lambda\,u - v =  f_3, \\
i\,\lambda\,\rho_u\,v - a_3\,u_{xx} - a_2\,z_{xx} = \rho_u\,f_4, \\
i\,\lambda\,\varphi + (y^{2} + \kappa)\,\varphi - w(w)\,\mu(y) = f_5,
\end{array}
\right.
\end{eqnarray}
with the boundary conditions
\begin{eqnarray}\label{513}
\begin{array}{l}
z(0)=z(L)=u(0)=u(L)=0.
\end{array}
\end{eqnarray}
Suppose that we have found $z$ and $u$ with the appropriated regularity. Therefore, from (\ref{512})$_{1,\,3}$ we have
\begin{eqnarray*}
\left\lbrace
\begin{array}{l}
w = i\,\lambda\,z - f_{1}, \\
v = i\,\lambda\,u - f_{3},
\end{array}
\right.
\end{eqnarray*}
it is clear that $z, u \in H_{0}^{1}(0, L)$. From (\ref{512})$_{5}$ we have
\begin{eqnarray}
\label{514}
\varphi(y) = \frac{f_{5}(x, y)\ +\ w(x)\,\mu(y)}{y^{2}\ +\ \kappa\ +\ i\,\lambda}.
\end{eqnarray}
On the other hand, replacing (\ref{512})$_{1,\,3}$ into (\ref{512})$_{2,\,4}$ respectively, we get
\begin{eqnarray}\label{516}
\left\{\matrix{
-\lambda^{2}\,\rho_z\,z - a_1\,z_{xx} - a_2\,u_{xx}+i\gamma\lambda(i\lambda+\kappa)^{\alpha-1}z(x)
= \rho_z\,f_{2} + i\rho_z\,\lambda\,f_{1}\hfill & \cr
\hfill+\gamma (i\lambda+\kappa)^{\alpha-1}f_1-\zeta \int_{-\infty}^{+\infty}\frac{\mu(y)f_{5}(x, y)}{y^2+i\lambda+\kappa}dy,  & \cr
-\lambda^{2}\,\rho_u\,u - a_3\,u_{xx} - a_2\,z_{xx} = \rho_u\,f_{4} + i\,\lambda\,f_{3}.\hfill & \cr}
\right.
\end{eqnarray}
Solving system (\ref{516}) is equivalent to find $z,\ u \in H^{2}(0, L)\cap H_{0}^{1}(0, L)$ such that
\begin{eqnarray}\label{517}
\left\{\matrix{(-\rho_z\lambda^{2}+\gamma i\lambda(i\lambda+\kappa)^{\alpha-1}) \Int_{0}^{L} z\ {\overline{\tilde z}}\ dx
 + a_1 \Int_{0}^{L}z_{x}\,{\overline{\tilde z}}_{x}\ dx+ a_2 \Int_{0}^{L}u_x\,{\overline{\tilde z}}_x\ dx\hfill &\cr
= \Int_{0}^{L}\left( \rho_z\,f_{2} + (\rho_z i\lambda+\gamma (i\lambda+\kappa)^{\alpha-1})\,\,f_{1} \right) {\overline{\tilde z}}\ dx
-\zeta \int_{0}^{L}\tilde{z}\int_{-\infty}^{+\infty}\frac{\mu(y)f_{5}(x, y)}{y^2+i\lambda+\kappa}dy\, dx,\hfill &\cr
-\lambda^{2}\,\rho_u \Int_{0}^{L} u\ {\overline{\tilde{u}}}\ dx
+ a_3 \Int_{0}^{L}u_{x}\,{\overline{\tilde{u}}}_{x}\ dx
+ a_2 \Int_{0}^{L}{\overline{\tilde{u}}}_x\,z_x\ dx
=\Int_{0}^{L} \left(\rho_u\,f_{4} + i\lambda\,\rho_u\,f_{3} \right){\overline{\tilde{u}}}\ dx\hfill &\cr}\right.
\end{eqnarray}
for all ${\cal V} = (\tilde{z},\ \tilde{u}) \in H_{0}^{1}(0, L)\times H_{0}^{1}(0, L)$.
The system (\ref{517}) is equivalent to the problem
\begin{eqnarray}
\label{520}
-\ \langle{{\cal L}}_{\lambda}\U,\,{\cal V}\rangle_{\left[H_{0}^{1}(0, L)\right]^{2}}
+ \langle \U,\, {\cal V}\rangle_{\left[H_{0}^{1}(0, L)\right]^{2}} = \Phi({\cal V}),
\end{eqnarray}
where
\begin{eqnarray*}
\matrix{\langle{{\cal L}}_{\lambda}\U,\,{\cal V}\rangle_{\left[H_{0}^{1}(0, L)\right]^{2}}
= \lambda^{2} \Int_{0}^{L}\left(\rho_z\,z\,{\overline{\tilde z}} + \rho_u\,u\,{\overline{\tilde u}} \right)dx}
\end{eqnarray*}
and
\begin{eqnarray*}
\matrix{\langle \U,\,{\cal V}\rangle_{\left[H_{0}^{1}(0, L)\right]^{2}}
= a_1 \Int_{0}^{L}z_{x}\,{\overline{\tilde z}}_{x}\ dx + a_3 \Int_{0}^{L}u_{x}\,{\overline{\tilde u}}_{x}\ dx + a_2 \Int_{0}^{L}\left(u_x\,{\overline{\tilde z}}_x + z_x\,{\overline{\tilde u}}_x \right)dx\hfill  \cr
+\gamma i\lambda(i\lambda+\kappa)^{\alpha-1}) \Int_{0}^{L} z\ {\overline{\tilde z}}\ dx. \cr}
\end{eqnarray*}
One can easily see that ${\cal L}_{\lambda},  \langle\ .,\, .\ \rangle_{\left[H_{0}^{1}(0, L)\right]^{2}}$ and
$\Phi$ are bounded. Furthermore
$$
\matrix{\Re\langle\U,\,\U\rangle_{\left[H_{0}^{1}(0, L)\right]^{2}}&=&
a_1 \Int_{0}^{L}|z_{x}|^2\, dx + a_3 \Int_{0}^{L}|u_{x}|^2\, dx + a_2 \Re\Int_{0}^{L}\left(u_x\,{\overline z}_x
 + z_x\,{\overline u}_x \right)dx\hfill  \cr
&&+\gamma \Re(i\lambda(i\lambda+\kappa)^{\alpha-1})) \Int_{0}^{L} |z|^2\ dx \cr
&\geq & a_1 \Int_{0}^{L}|z_{x}|^2\, dx + a_3 \Int_{0}^{L}|u_{x}|^2\, dx + a_2 \Re\Int_{0}^{L}\left(u_x\,{\overline z}_x
+ z_x\,{\overline u}_x \right)dx  \cr
&= & (a_1-\frac{a_2^2}{a_3}) \Int_{0}^{L}|z_{x}|^2\, dx + \Int_{0}^{L}\left|\frac{a_2}{\sqrt{a_3}}z_{x}+\sqrt{a_3}u_{x}\right|^2\, dx,\hfill \cr}
$$
where we have used the fact that
$$
\gamma\lambda\Re \left(i(i\lambda+\omega)^{\alpha-1}\right)=\zeta\lambda^2\Int_{-\infty}^{+\infty}\Frac{\mu(y)^2}{\lambda^2+(\omega+y^2)^2}\, dy> 0.
$$
Thus $\langle\cdot, \cdot\,\rangle_{\left[H_{0}^{1}(0, L)\right]^{2}}$ is coercive. Consequently, by Fredholm alternative, proving
the existence of $U$ solution of (\ref{520}) reduces to proving that (\ref{520}) with $\Phi\equiv 0$ has a notrivial solution. Indeed
if there exists $U\not=0$, such that
\begin{equation}
L_{\lambda}(U, V)+a_{(H_{*}^{1}(0 , L))^2}(U, V)=0\quad \forall V\in H_{0}^{1}(0 , L)\times H_{0}^{1}(0, L).
\label{eaz}
\end{equation}
In particular for $V=U$, it follows that
$$
\matrix{\lambda^{2} \Int_{0}^{L}\left(\rho_z\,|z|^{2} + \rho_u\,|u|^{2}\right)dx+
a_1 \Int_{0}^{L}|z_{x}|^2\, dx + a_3 \Int_{0}^{L}|u_{x}|^2\, dx + a_2 \Int_{0}^{L}\left(u_x\,{\overline z}_x
+ z_x\,{\overline u}_x \right)dx\hfill & \cr
+\gamma (i\lambda(i\lambda+\kappa)^{\alpha-1})) \Int_{0}^{L} |z|^2\ dx=0. & \cr}
$$
Hence, we have
\begin{equation}
u=v=0.
\label{e22xx1}
\end{equation}
We deduce that $U= 0$. Hence $i\lambda- {\cal A}$ is surjective for all $\lambda\in \R^*$.\\
\begin{lemma}
If $\kappa \neq 0,$ we have that $\ 0\in\varrho({{\cal A}}).$
\label{lemma2}
\end{lemma}

{\bf Proof.}
Now, given $F=(f_{1},\,f_{2},\,f_{3},\,f_{4},\,f_{5})^{T} \in {{\cal H}}$, we must show there exists a unique $U=(z,\,w,\,u,$ $\,v,\,\varphi)^{T}$ in ${{\cal D}}({{\cal A}}),$ such that $-{{\cal A}}U=F$, namely,
\begin{eqnarray}\label{526}
\left\lbrace
\begin{array}{l}
- w = f_{1}, \\
- a_1\,z_{xx} -a_2\,u_{xx}+\zeta \int_{-\infty}^{+\infty}\mu(y)\varphi(x,y)dy = \rho_z\,f_{2}, \\
- v = f_{3}, \\
- a_3\,u_{xx} -a_2\,z_{xx} = \rho_u\,f_{4}, \\
(y^{2} + \kappa)\,\varphi - w(x)\,\mu(y) = f_{5},
\end{array}
\right.
\end{eqnarray}
with the boundary conditions
$$
z(0)=z(L)=u(0)=u(L)=0.
$$
Using $(\ref{526})_2$ and $(\ref{526})_4$ with $\lambda = 0$ it follows that
\begin{eqnarray}\label{528}
\matrix{a_1 \Int_{0}^{L}z_{x}\,{\overline{\tilde z}}_{x}\ dx + a_2 \Int_{0}^{L}u_x\,{\overline{\tilde z}}_x\ dx
&= & \, \rho_z \Int_{0}^{L}\,f_{2}\,{\overline{\tilde z}}\ dx
+\gamma \kappa^{\alpha-1}\Int_{0}^{L}f_1 {\overline{\tilde z}}\ dx\hfill \cr
&&-\zeta\Int_{0}^{L}{\overline{\tilde z}}\Int_{-\infty}^{+\infty}\frac{\mu(y)\,f_{5}(x, y)}{y^{2} + \kappa}\ dy\, \ dx&\cr}
\end{eqnarray}
and
\begin{eqnarray}\label{5280}
a_3 \Int_{0}^{L}u_{x}\,{\overline{\tilde u}}_{x}\ dx + a_2 \Int_{0}^{L}z_x\,{\overline{\tilde u}}_x\ dx
 = & \, \rho_u \Int_{0}^{L}f_{4}\,{\overline{\tilde u}}\ dx
\end{eqnarray}
for all $\tilde{z}, \tilde{u} \in H_{0}^{1}(0, L)$. The system (\ref{528})--(\ref{5280}) is equivalent to the problem
\begin{eqnarray}\label{529}
a_{\kappa}\left((z,\,u),\ (\tilde{z},\,\tilde{u})\right) = {{\cal L}}_{\kappa}(\tilde{z},\,\tilde{u}),
\end{eqnarray}
where the continuous sesquilinear form  and coercive $a_{\kappa}: \left[H_{0}^{1}(0, L)\times
H_{0}^{1}(0, L)\right]^{2}\rightarrow\C$ and the continuous antilinear form
${{\cal L}}_{\kappa}: \left[H_{0}^{1}(0, L)\right]^{2}\rightarrow\C$ are defined by
\begin{equation} \label{530}
a_{\kappa}\big((z,\,u),\ (\tilde{z},\,\tilde{u})\big) = a_1 \Int_{0}^{L}z_{x}\,{\overline{\tilde z}}_{x}\ dx + a_3 \Int_{0}^{L}u_{x}\,{\overline{\tilde u}}_{x}\ dx
+ a_2 \Int_{0}^{L} \left( u_x\,{\overline{\tilde z}}_x + z_x\,{\overline{\tilde u}}_x \right)\ dx
\end{equation}
and
\begin{eqnarray}\label{531}
\matrix{{{\cal L}_{\kappa}}(\tilde{z},\,\tilde{u}) = & \rho_z \Int_{0}^{L}f_{2}\,{\overline{\tilde z}}\ dx
+ \rho_u \Int_{0}^{L}f_{4}\,{\overline{\tilde u}}\ dx
+\gamma \kappa^{\alpha-1}\Int_{0}^{L}f_1 {\overline{\tilde z}}\ dx\hfill \cr
& -\zeta\Int_{0}^{L}{\overline{\tilde z}}\Int_{-\infty}^{+\infty}\frac{\mu(y)\,f_{5}(x, y)}{y^{2} + \kappa}\ dy\, \ dx.
\hfill \cr}
\end{eqnarray}
Applying the Lax-Milgram Theorem, we have that, for all $(\tilde{z},\,\tilde{u})\in H_{0}^{1}(0, L)\times H_{0}^{1}(0, L)$
the problem (\ref{529}) admits a unique solution $(z,\,u)\in H_{0}^{1}(0, L)\times H_{0}^{1}(0, L)$.
Using elliptic regularity, it follows from (\ref{528})-(\ref{5280}) that $(z,\,u) \in H^{2}(0, L)\times H^{2}(0, L)$.
Therefore, the operator ${{\cal A}}$ is surjective.\\

\subsection{Polynomial stability and optimality (for $\kappa\not= 0$)}
To achieve proof of polynomially stable of C$_{0}$-semigroup $e^{t{{\cal A}}}$ we use the following result.
\begin{theorem}\cite{bo} \label{bo}
Let ${{\cal S}}(t) = e^{{{\cal A}}\,t}$ be a C$_{0}$-semigroup of contractions on Hilbert space ${{\cal H}}.$ If
\begin{eqnarray*}
i\,\R\subseteq\varrho({{\cal A}})\quad\mbox{and}\quad \sup_{|\beta|\geq 1}\frac{1}{\beta^{\ell}}\,\|(i\,\beta\,I - {{\cal A}})^{-1}\|_{{{\cal L}}({{\cal H}})} < M,
\end{eqnarray*}
for some $\ell$, then there exist $c$ such that
\begin{eqnarray*}
\|e^{{{\cal A}}\,t}U_{0}\|^{2} \leq \frac{c}{t^{\frac{2}{\ell}}}\,\|U_{0}\|_{{{\cal D}}({{\cal A}})}^{2}.
\end{eqnarray*}
\end{theorem}

Now, we will prove the second main theorem of this section.
\begin{theorem}\label{theorem507}
The semigroup ${{\cal S}}_{{{\cal A}}}(t)_{t\geq 0}$ is polynomially stable and
\begin{eqnarray}\label{536}
{\cal E}(t)=\|{{\cal S}}_{{{\cal A}}}(t)U_{0}\|_{{{\cal H}}}^2 \leq \frac{1}{t^{2/\,(1 - \alpha)}}\,\|U_{0}\|_{{{\cal D}}({{\cal A}})}^2.
\end{eqnarray}
Moreover, the rate of energy decay $t^{-2/(1-\alpha)}$ is optimal for general initial data in $D({\cal A})$.

\end{theorem}

{\bf Proof.}
In section 3, we have proved that the first condition in Theorem \ref{bo} is satisfied. Now, we need to show that
\begin{equation}
\Sup_{|\lambda|\geq 1}\Frac{1}{\lambda^{l}}\|(i\lambda I-{\cal A})^{-1}\|_{{\cal H}}< \infty,
\label{v35}
\end{equation}
where $l=1-\alpha$.
We establish (\ref{v35}) by contradiction. So, if (\ref{v35}) is false, then there exist sequences
$(\lambda_n)_n\subset \R$ and $U_n= (u_n, {\tilde u}_n, v_n, {\tilde v}_n, \varphi_n)\in  D({\cal A})$
satisfying
\begin{equation}
\|U_n\|_{{\cal H}}=1\quad \forall n\geq 0,
\label{v37}
\end{equation}
\begin{equation}
\Lim_{n\rightarrow\infty}|\lambda_n|=+\infty
\label{v38}
\end{equation}
and
\begin{equation}
\Lim_{n\rightarrow\infty}\lambda_n^{l}\|(i\lambda_n I-{\cal A})U_n\|\rightarrow 0,
\label{v39}
\end{equation}
which implies that
\begin{eqnarray}\label{537}
\left\lbrace
\begin{array}{l}
\lambda^{l}(i\,\lambda\,z - w)= g_{1}\rightarrow 0 \hbox{ in } H_{0}^1(0, L), \\
\lambda^{l}(i\,\lambda\,\rho_z\,w - a_1\,z_{xx} - a_2\,u_{xx}
+\zeta \int_{-\infty}^{+\infty}\mu(y)\varphi(x,y)dy) = \rho_z\,g_{2}\rightarrow 0 \hbox{ in } L^2(0, L),  \\
\lambda^{l}(i\,\lambda\,u - v) = g_{3}\rightarrow 0 \hbox{ in } H_{0}^1(0, L), \\
\lambda^{l}(i\,\lambda\,\rho_u\,v - a_3\,u_{xx} - a_2\,z_{xx}) = \rho_u\,g_{4}\rightarrow 0 \hbox{ in } L^2(0, L),  \\
\lambda^{l}((y^{2} + \kappa + i\,\lambda)\,\varphi - w(x)\,\mu(y)) = g_{5}\rightarrow 0 \hbox{ in } L^2(-\infty, +\infty),
\end{array}
\right.
\end{eqnarray}
For simplification, we denote $\lambda_n$ by $\lambda, U_n= (z_n, w_n, u_n, v_n, \varphi_n)$ by
$U= (z, w, u, v, \varphi)$ and
$G_n=(g_{1n}, g_{2n}, g_{3n}, g_{4n}, g_{5n})=\lambda_n^{l}(i\lambda_n I-{\cal A})U_n$ by
$G=(g_{1}, g_{2}, g_{3}, g_{4}, g_{5})$.
We will prove that\\ $\|U\|_{\cal H}=o(1)$ as a contradiction with (\ref{v37}).
Our proof is divided into several steps.

\noindent
$\bullet${\bf Step 1}
Taking the inner product of $\lambda^{l}(i\lambda I-{\cal A})U$ with $U$, we get
\begin{equation}
i\lambda \|U\|_{{\cal H}}^{2}-({\cal A}U, U)_{{\cal H}}=\Frac{o(1)}{\lambda^{l}}.
\label{v2}
\end{equation}
Using (\ref{406}), we get
\begin{equation}
\zeta\int_{0}^{L}\int_{-\infty}^{+\infty}(y^{2}+\omega)|\varphi(y)|^{2}\, dy\, dx
=-\Re({\cal A}U, U)=\Frac{o(1)}{\lambda^{l}}.
\label{v3}
\end{equation}
From (\ref{537})$_{5}$ we obtain
\begin{eqnarray}\label{541}
w(x)\,\mu(y) = (y^{2} + \kappa + i\,\lambda)\,\varphi - f_{5}(y).
\end{eqnarray}
By multiplying $(\ref{541})$ by $(i\lambda+y^{2}+\kappa)^{-2}|y|$, we get
\begin{equation}
(i\lambda+y^{2}+\kappa)^{-2}w(x)\mu(y)|y|=(i\lambda+y^{2}+\kappa)^{-1}|y|\varphi-(i\lambda+y^{2}+\kappa)^{-2}|y|f_5(x, y).
\label{e38kk}
\end{equation}
Hence, by taking absolute values of both sides of (\ref{e38kk}), integrating over the interval $]-\infty, +\infty[$ with
respect to the variable $y$ and applying Cauchy-Schwartz inequality, we obtain
\begin{equation}
{\cal S}|w(x)|\leq  \sqrt{2}{\cal U} \left(\Int_{-\infty}^{+\infty}y^{2}|\varphi|^{2}\, dy\right)^{\frac{1}{2}}
+ 2 {\cal V}\left(\Int_{-\infty}^{+\infty}|f_5(x, y)|^{2}\, dy\right)^{\frac{1}{2}},
\label{e39kk}
\end{equation}
where
$$
{\cal S}=\left|\Int_{-\infty}^{+\infty}(i\lambda+y^{2}+\kappa)^{-2}|y|\mu(y)\, dy\right|
=\frac{|1-2\alpha|}{4}\Frac{\pi}{|\sin\frac{(2\alpha+3)}{4}\pi|}|i\lambda+\kappa|^{\frac{(2\alpha-5)}{4}},
$$
$$
{\cal U}=\left(\Int_{-\infty}^{+\infty}(|\lambda|+y^{2}+\kappa)^{-2}\, dy\right)^{\frac{1}{2}}=(\frac{\pi}{2})^{1/2}||\lambda|+\kappa|^{-\frac{3}{4}},
$$
$$
{\cal V}=\left(\Int_{-\infty}^{+\infty}(|\lambda|+y^{2}+\kappa)^{-4}|y|^{2}\, dy\right)^{\frac{1}{2}}
=\left(\Frac{\pi}{16}||\lambda|+\kappa|^{-\frac{5}{2}}\right)^{1/2}.
$$
Thus, by using the inequality $2PQ \leq P^2 + Q^2, P \geq 0, Q \geq 0$, again, we get
\begin{equation}
\ \ \ \ \ \ \ \ {\cal S}^2\|w(x)\|_{L^2(0, L)}^2\leq  2 {\cal U}^2 \left(\Int_{0}^{L}\Int_{-\infty}^{+\infty}(y^{2}+\kappa)|\varphi|^{2}\, dy\, dx\right)
+ 4{\cal V}^2\left(\Int_{0}^{L}\Int_{-\infty}^{+\infty}|f_5(x, y)|^{2}\, dy\, dx\right).
\label{e339kk}
\end{equation}
We deduce that
\begin{equation}
\|w(x)\|_{L^2(0, L)}^2=\Frac{o(1)}{\lambda^{l-(1-\alpha)}}+\Frac{o(1)}{\lambda^{2l+\alpha}}.
\label{e41nkk}
\end{equation}
Then
\begin{equation}
\|w(x)\|_{L^2(0, L)}=\Frac{o(1)}{\lambda^{\frac{l-(1-\alpha)}{2}}}.
\label{mn20}
\end{equation}
By eliminating $w$ and $v$ from system (\ref{537}) we obtain
\begin{equation}
-\lambda^2\,\rho_z z - a_1\,z_{xx} - a_2\,u_{xx}
+\zeta \int_{-\infty}^{+\infty}\mu(y)\varphi(x,y)dy= \rho_z\,f_{2}+i\lambda\rho_zf_1  \hbox{ in } L^2(0, L),
\label{v7}
\end{equation}
\begin{equation}
-\lambda^2\rho_u u - a_3\,u_{xx} - a_2\,z_{xx}=\rho_u\,f_{4}+i\lambda\rho_u f_3   \hbox{ in } L^2(0, L),
\label{v8}
\end{equation}
where
\begin{equation}
\left\{\matrix{\|f\|_{L^2(0, 1)}=\left\|\Frac{\rho_z g_2+i\lambda\rho_z g_1}{\lambda^l}\right\|_{L^2(0, 1)}=\Frac{o(1)}{\lambda^{l-1}},\hfill & \cr
\|g\|_{L^2(0, 1)}=\left\|\Frac{\rho_u g_4+i\lambda\rho_z g_3}{\lambda^l}\right\|_{L^2(0, 1)}=\Frac{o(1)}{\lambda^{l-1}}.\hfill & \cr}\right.
\label{v9}
\end{equation}
Note that thanks to $a_1+a_2\nu> 0$ (where $\nu=-a_2/a_3$), (\ref{v7}) and (\ref{v8}) can be transformed into the following
\begin{equation}
-\lambda^2\,\rho_z z - (a_1+a_2\nu)\,z_{xx} +a_2(\nu\,z_{xx}-\,u_{xx})
+\zeta \int_{-\infty}^{+\infty}\mu(y)\varphi(x,y)dy= \rho_z\,f_{2}+i\lambda\rho_zf_1,
\label{vv7}
\end{equation}
\begin{equation}
-\lambda^2\nu\rho_u u -a_2(\nu\,z_{xx}-\,u_{xx}) =\nu(\rho_u\,f_{4}+i\lambda\rho_u f_3)   \hbox{ in } L^2(0, L),
\label{vv8}
\end{equation}
Thus, adding (\ref{vv7}) and (\ref{vv7}) so as to eliminate the common item $a_2(\nu\,z_{xx}-\,u_{xx})$
\begin{equation}
\matrix{\lambda^2\,\rho_z z+\lambda^2\,\rho_u \nu u + (a_1+a_2\nu)\,z_{xx}
-\zeta \int_{-\infty}^{+\infty}\mu(y)\varphi(x,y)dy= -(\rho_z\,f_{2}+i\lambda\rho_zf_1)\hfill \cr
\hspace{3cm}-\nu(\rho_u\,f_{4}+i\lambda\rho_u f_3)  \hbox{ in } L^2(0, L),\cr}
\label{vv77}
\end{equation}
Then, taking the $L^2$-inner product of (\ref{vv77}) with $z$, we obtain
\begin{equation}
\matrix{\lambda^2\,\rho_z \|z\|_{L^2(0, L)}^2+\lambda^2\,\rho_u \Int_{0}^{L}\nu u {\overline z}\, dx - (a_1+a_2\nu)\|z_{x}\|_{L^2(0, L)}^2
-\zeta \Int_{0}^{L}{\overline z}\int_{-\infty}^{+\infty}\mu(y)\varphi(x,y)dy\, dx\hfill &\cr
= -\Int_{0}^{L}{\overline z}[(\rho_z\,f_{2}+i\lambda\rho_zf_1)+\nu(\rho_u\,f_{4}+i\lambda\rho_u f_3)]\, dx  \hbox{ in } L^2(0, L),\hfill &\cr}
\label{vv78}
\end{equation}
By (\ref{vv77}), (\ref{vv77}) and Cauchy-Schwartz inequality, we get that the second term in (\ref{vv78}) satisfies
\begin{equation}
|\lambda^2\,\rho_u \Int_{0}^{L}\nu u {\overline z}\, dx|^2\sim |\rho_u \nu \Int_{0}^{L}v {\overline w}\, dx|^2
\leq (\rho_u \nu)^2\Int_{0}^{L}|v|^2\, dx\Int_{0}^{L}|w|^2\, dx=\Frac{o(1)}{\lambda^{l-(1-\alpha)}}.
\label{vv79}
\end{equation}
Here we have used (\ref{mn20}) and the boundedness of $\|v\|_{L^2(0, L)}$.
Also
\begin{equation}
\Int_{0}^{L}{\overline z}[(\rho_z\,f_{2}+i\lambda\rho_zf_1)+\nu(\rho_u\,f_{4}+i\lambda\rho_u f_3)]\, dx
=\Frac{o(1)}{\lambda^{l+\frac{l-(1-\alpha)}{2}}}.
\label{vv80}
\end{equation}
We can estimate
$$
\matrix{\left|\Int_{0}^{L}\bar{z}(\Int_{-\infty}^{+\infty}\mu(y)\varphi(x,y)\, dy)\, dx\right|\hfill &\cr
\leq \|z\|_{L^2(0, L)}\left(\Int_{-\infty}^{+\infty}\Frac{\mu^2(y)}{y^2+\kappa}\, dy\right)^{\frac{1}{2}}
\left(\Int_{0}^{L}\Int_{-\infty}^{+\infty}(y^2+\kappa)|\varphi(x,y)|^2\, dy\, dx\right)^{\frac{1}{2}}
=\Frac{o(1)}{\lambda^{\frac{l}{2}}}.\hfill &\cr}
$$
This together with (\ref{vv78}), (\ref{vv79}) and (\ref{vv80}) shows that
\begin{equation}
\lambda^2\,\rho_z \|z\|_{L^2(0, L)}^2 - (a_1+a_2\nu)\|z_{x}\|_{L^2(0, L)}^2=\Frac{o(1)}{\lambda^{\frac{l-(1-\alpha)}{2}}}.
\label{vv81}
\end{equation}
which along with $(\ref{537})_1$ and (\ref{mn20}) leads to
\begin{equation}
\|z_{x}\|_{L^2(0, L)}^2=\Frac{o(1)}{\lambda^{\frac{l-(1-\alpha)}{2}}}.
\label{vv82}
\end{equation}
Taking the $L^2$-inner product of (\ref{vv77}) with $u$, we obtain
\begin{equation}
\matrix{\lambda^2\,\rho_z \Int_{0}^{L}\nu {\overline u} z\, dx+\lambda^2\,\rho_u \nu \|u\|_{L^2(0, L)}^2
-(a_1+a_2\nu)\Int_{0}^{L}z_{x}{\overline u}_x\, dx
-\zeta \Int_{0}^{L}{\overline u}\int_{-\infty}^{+\infty}\mu(y)\varphi(x,y)dy\, dx\hfill &\cr
= -\Int_{0}^{L}{\overline u}[(\rho_z\,f_{2}+i\lambda\rho_zf_1)+\nu(\rho_u\,f_{4}+i\lambda\rho_u f_3)]\, dx  \hbox{ in } L^2(0, L).\hfill &\cr}
\label{vv83}
\end{equation}
On the one hand, we have
$$
\lambda^2\,\rho_z \Int_{0}^{L}\nu {\overline u} z\, dx=\Frac{o(1)}{\lambda^{\frac{l-(1-\alpha)}{2}}}.
$$
On the other hand, using H\"{o}lder inequality again, along with (\ref{vv82}) and the boundedness
of $\|u_x\|_{L^2(0, L)}$, we obtain that
$$
\Int_{0}^{L}z_{x}{\overline u}_x\, dx=\Frac{o(1)}{\lambda^{\frac{l-(1-\alpha)}{4}}}.
$$
Hence, by (\ref{vv83}), we get
$$
\lambda^2\,\rho_u \nu \|u\|_{L^2(0, L)}^2=\Frac{o(1)}{\lambda^{\frac{l-(1-\alpha)}{4}}}.
$$
which along with $(\ref{537})_3$ leads to
\begin{equation}
\|v\|_{L^2(0, L)}=\Frac{o(1)}{\lambda^{\frac{l-(1-\alpha)}{8}}}.
\label{vv84}
\end{equation}
Taking the $L^2$-inner product of (\ref{v7}) and (\ref{v8}) with $z$ and $u$, respectively,
then integrating by parts yields that
\begin{equation}
\matrix{-\lambda^2\,\rho_z \|z\|_{L^2(0, L)}^2+a_1\|z_{x}\|_{L^2(0, L)}^2+a_2 \Int_{0}^{L} u_x {\overline z}_x\, dx
+\zeta \Int_{0}^{L}{\overline z}\int_{-\infty}^{+\infty}\mu(y)\varphi(x,y)dy\, dx\hfill &\cr
= \Int_{0}^{L}{\overline z}(\rho_z\,f_{2}+i\lambda\rho_zf_1)\, dx\ \   \hbox{ in } L^2(0, L)\hfill &\cr}
\label{vv85}
\end{equation}
and
\begin{equation}
\matrix{-\lambda^2\,\rho_u \|u\|_{L^2(0, L)}^2++a_3\|u_{x}\|_{L^2(0, L)}^2+a_2 \Int_{0}^{L} z_x {\overline u}_x\, dx
\hfill &\cr
= \Int_{0}^{L}{\overline u}(\rho_u\,f_{4}+i\lambda\rho_u f_3)\, dx\ \   \hbox{ in } L^2(0, L). &\cr}
\label{vv86}
\end{equation}
Adding (\ref{vv85}) and (\ref{vv86}), and by virtue of  $a_1+\nu a_2>0$, we obtain
\begin{equation}
\matrix{-\lambda^2\,\rho_z \|z\|_{L^2(0, L)}^2
-\lambda^2\,\rho_u \|u\|_{L^2(0, L)}^2+(a_1+\nu a_2)\|z_{x}\|_{L^2(0, L)}^2\hfill &\cr
+\Int_{0}^{L}|\frac{a_2}{\sqrt{a_3}}z_x+\sqrt{a_3}u_x|^2\, dx
+\zeta \Int_{0}^{L}{\overline z}\int_{-\infty}^{+\infty}\mu(y)\varphi(x,y)dy\, dx
\hfill &\cr
= \Int_{0}^{L}[{\overline z}(\rho_z\,f_{2}+i\lambda\rho_zf_1)+{\overline u}(\rho_u\,f_{4}+i\lambda\rho_u f_3)]\, dx
\ \ \hbox{ in } L^2(0, L).\hfill &\cr}
\label{vv87}
\end{equation}
Substituting (\ref{mn20}) and (\ref{vv84}) into (\ref{vv87})
\begin{eqnarray}\label{1}
\Int_{0}^{L}\left(\rho_z\,|w|^{2} + \rho_u\,|v|^{2} + (a_1-\frac{a_2^2}{a_3})\,|z_{x}|^{2}
+ \left|\frac{a_2}{\sqrt{a_3}}z_{x}+\sqrt{a_3} u_x\right|^2\right)dx=o(1).
\end{eqnarray}
Moreover, from (\ref{v3}) we obtain
\begin{eqnarray}
\displaystyle{\int_{\R}}|\varphi|^{2}\ dy \leq \Frac{1}{\kappa}\int_{\R}\left(y^{2} + \kappa\right)|\varphi|^{2}\ dy
=\Frac{o(1)}{\lambda^{l}}.
\label{552}
\end{eqnarray}
Thus, we conclude that
$$
\|U\|_{{\cal H}}=o(1).
$$
Besides, we prove that the decay rate is optimal. Indeed, the decay rate is consistent
with the asymptotic expansion of eigenvalues which shows a behavior of the real part
like $\mu_n^{-(1-\alpha)}$.

\section*{Declarations}
{\bf Conflict of interest} The authors declare that they have no conflict of interest.

\end{sloppypar}


\begin{thebibliography}{20}
\bibitem{acbe} Z. Achouri, N. Amroun and A. Benaissa, {\it The Euler-Bernoulli beam equation with
boundary dissipation of fractional derivative type\/,} Mathematical Methods in the Applied Sciences
{\bf 40}, (2017)-11, 3837-3854.
\bibitem{Santos} R.G.C. Almeida, M.L. Santos, {\it Lack of exponential decay of a coupled system of wave equations
with memory\/,} Nonlinear Anal-Real.  {\bf 12} (2011) 1023-1032.
\bibitem{apala} Apalara, T.A. {\em General stability result of swelling porous elastic soils with a
viscoelastic damping\/,} Z. Angew. Math. Phys. {\bf 71}-200, (2020).
\bibitem{amma} K. Ammari, H. Fathi and L. Robbiano, {\em Fractional-feedback stabilization for a class of evolution
systems\/,} Journal of Differential Equations, {\bf 268} (2020)-1, 5751-5791.
\bibitem{ramos1} J. A. Anderson, Ramos, Cledson S. L. Gonçalves, N. Corrêa, S. Silv\'erio,
{\em Exponential stability and numerical treatment for piezoelectric beams with magnetic effect\/,}
ESAIM Math. Model. Numer. Anal. {\bf 52} (2018)-1, 255-274.
\bibitem{arendt} W. Arendt, C. J. K. Batty, {\it Tauberian theorems and stability of one-parameter semigroups\/,}
Trans. Amer. Math. Soc.  {\bf 306} (1988)-2, 837-852.
\bibitem{menno} T. Bentrcia, A. Mennouni, {\em On the asymptotic stability of a Bresse system with two
fractional damping terms: Theoretical and numerical analysis\/,} Discrete and Continuous
Dynamical Systems - Series B. {\bf 28} (2023), 580-622.
\bibitem{bo} A. Borichev, Y. Tomilov, {\it Optimal polynomial decay of functions and operator semigroups,}
Math. Ann. Vol. {\bf 347} (2010)-2, 455-478.
\bibitem{boy} L. Boyadjiev, O. Kamenov,  S.L. Kalla, {\it On the Lauwerier formulation of the temperature field problems
in oil strata\/,} International J. Math. Math. Sci. {\bf 10} (2005) 1577-1588.
\bibitem{CM} J. Choi, R.  Maccamy, {\it Fractional order Volterra equations with applications to elasticity\/,}
J. Math. Anal. Appl. {\bf 139} (1989) 448-464.
\bibitem{Cordeiro_Lobato_Raposo} S. Cordeiro, R. F. C. Lobato, C. A. Raposo, {\it Optimal polynomial decay for a coupled
system of wave with past history\/,} Open J. Math. Anal. {\bf 4} (2020) 49-59.
\bibitem{iesa} D. Iesan, {\em On the theory of mixtures of thermoelastic solids\/,} J. Thermal Stress. {\bf 14}
(1991)-4, 389–408.
\bibitem{KST} A. Kilbas, H. Srivastava, J. Trujillo,
{\it Theory and Applications of Fractional Differential Equations\/,} North-Holland Mathematics Studies, 204.
Elsevier Science B.V., Amsterdam, (2006).
\bibitem{111} A. Kilbas, and J. Trujillo, {\it Differential equation  of fractional order: methods, results and Problems\/,} Appl. Anal. Vol.I 78(2) (2002) 435-493.
\bibitem{222} A. Kilbas, and J. Trujillo, {\it Differential equation  of fractional order: methods, results and problems\/,}
Appl. Anal. Vol. Vol.II,  {\bf 81} (2001)-(1-2) 153-192.
\bibitem{Renato} R.F.C. Lobato, S.M.S. Cordeiro, M.L. Santos, D.S. Almeida Junior, {\it Optimal polynomial decay to
coupled wave equations and its numerical properties\/,} J. Appl. Math.  (2014).  Art. Id 897080
\bibitem{lyub} I. Lyubich Yu \& V.Q. Phóng, {\em Asymptotic stability of linear differential equations in Banach
spaces\/,} {Studia Mathematica\/,} {\bf 88} (1988)-(1), 37-42.
\bibitem{Najafi} M. Najafi, {\it Study of exponential stability of coupled wave systems via distributed stabilizer\/,}
Int. J. Math. Math. Sci. {\bf 28} (2001) 479-491.
\bibitem{nonato}  C. Nonato, A. Benaissa, A. Ramos, C. Raposo, M. Freitas, {\em Porous elastic soils with fluid
saturation and dissipation of fractional derivative type\/,} Qual. Theory Dyn. Syst. {\bf 23}, Art. 79 (2024).
\bibitem{olive} W. Oliveira, S. Cordeiro, C. Raposo, O. Vera, {\em Asymptotic behavior for a porous-elastic system
with fractional derivative-type internal dissipation\/,} FCAA, DOI: 10.1007/s13540-024-00250-y.
\bibitem{Pereira_PJAA} D.C. Pereira, C.A. Raposo, C. Maranh\~ao, A. Cattai,
{\it Wave coupled system of the p-Laplacian type\/,} Poincar\'e J. Anal. Appl. {\bf 7} (2020) 185-195.
\bibitem{quint} R. Quintanilla, {\em Exponential stability for one-dimensional problem of swelling porous elastic
soils with fluid saturation\/,} J. Comput. Appl. Math., {\bf 145} (2002), 525-533.
\bibitem{santos} A. J. A. Ramos, M. M. Freitas, D. S. Almeida, Jr. A. S. No\'e and M. J. Dos Santos, {\em Stability
results for elastic porous media swelling with nonlinear damping\/,} J. Math. Phys., {\bf 61} (2020),101505, 10pp.
\bibitem{15} B. Mbodje, {\it Wave energy decay under fractional derivative controls\/,}
IMA. IMA J. Math. Control Inf. {\bf 23} (2006) 237-257.
\bibitem{SKM} S. Samko, A. Kilbas, O. Marichev, {\it Integral and Derivatives of Fractional Order\/,}
Gordon  Breach, New York, 1993.
\bibitem{guowa} J.M.Wang and B.-Z. Guo, {\em On the stability of swelling porous elastic soils with fluid saturation
by one internal damping\/,} IMA J. Appl. Math., {\bf 71} (2006), 565-582.















\end{thebibliography}
\end{document}